\documentclass[12pt,reqno]{amsart}

\usepackage{bm}
\usepackage{amssymb}
\usepackage{amsmath, amsfonts}
\usepackage{graphicx}
\usepackage{fullpage}
\newtheorem{theorem}{Theorem}
\newtheorem{remark}[theorem]{Remark}

\def\e{\varepsilon}
\def\N{{\mathcal N}}
\def\R{{\mathcal R}}
\def\U{{\mathcal U}}
\def\integer{{\mathbb Z}}
\def\real{{\mathbb R}}	
\def\torus{{\mathbb T}}

\newcommand{\Addresses}{{
\bigskip
\footnotesize

A.P.~Bustamante, \textsc{School of Mathematics, Georgia Institute of Technology, Atlanta, USA}\par\nopagebreak \textit{E-mail address}, \texttt{apb7@math.gatech.edu}

\medskip
C. Chandre, \textsc{CNRS, Aix Marseille Univ, I2M, 13009 Marseille, France}\par\nopagebreak \textit{E-mail address}, \texttt{chandre@math.cnrs.fr}

}}

\begin{document}

\title{ Numerical computation of Critical surfaces for the breakup of invariant tori in Hamiltonian systems }

\author{ Adri\'an P. Bustamante \and  Cristel Chandre } 


\begin{abstract}
We compute the critical surface for the existence of invariant tori of a family of Hamiltonian systems with two and three degrees of freedom. We use and compare two methods to compute the critical surfaces: renormalization-group transformations and conjugation in configuration space. We unveil the presence of cusps in the critical surface for the breakup of three-dimensional invariant tori, whereas the critical surface of two-dimensional invariant tori is expected to be smooth. 
\end{abstract}

\maketitle

\section{Introduction}

We consider Hamiltonian systems of the form
\begin{equation}
\label{eqn:Hamiltonian}
    H({\bf A}, {\bm \varphi})={\bm \omega}\cdot {\bf A} +\frac{1}{2}({\bm \Omega}\cdot {\bf A})^2+V({\bm \varphi}),
\end{equation}
where ${\bm \varphi}\in \mathbb{T}^d$ and ${\bf A}\in \mathbb{R}^d$. We are asking whether of not Hamiltonian~\eqref{eqn:Hamiltonian} has an invariant torus with frequency vector $\bm\omega$. If $V=0$, the equations of motion show that the Hamiltonian system~\eqref{eqn:Hamiltonian} possesses an invariant torus with frequency vector ${\bm\omega}$ at $z\equiv {\bm\Omega}\cdot {\bf A}=0$. For sufficiently small and regular $V$ and under suitable assumptions on the frequency vector ${\bm\omega}$, KAM theory ensures the persistence of this torus, which is slightly deformed by the perturbation.  
The proof revolves around a change of variables such that the flow associated with Hamiltonian~\eqref{eqn:Hamiltonian} is locally conjugated to a rotation by $\bm\omega$. In other words, there exists a local change of variables $({\bf A}, {\bm \varphi})\mapsto (\overline{\bf A}, {\bm \psi})$ such that $\dot{{\bm \psi}}={\bm\omega}$ and $\dot{\overline{\bf A}}={\bf 0}$ for $\overline{\bf A}={\bf 0}$.
Here we construct numerically this conjugation using two different methods: a renormalization-group transformation in the space of Hamiltonians, and a conjugation method in configuration space. 
In what follows, we consider potentials $V$ parameterized by two parameters. The set of these parameters for which the conjugation can be found (or equivalently, for which there is an invariant torus with frequency ${\bm\omega}$) is bounded by a surface, called critical surface. The comparison of the critical surface using the two methods of conjugation sheds light on the methods, and also highlights which features are artefacts of the method and which are actual features of the critical surface. 

We consider two- and three-degrees of freedom Hamiltonian systems, i.e., $d=2$ and $d=3$. In this article, we restrict frequency vectors to the ones where there exist a square integer matrix, $N$, with determinant $\pm 1$, such that $N{\bm\omega} = \theta_1{\bm\omega}$ with $|\theta_1| < 1$ and the remaining eigenvalues of $N$ such that $|\theta_i| > 1$, for $i>1$. These frequency vectors satisfy a Diophantine condition~\cite{Koch-99}. 

In Sec.~\ref{sec:renorm}, we detail the chosen renormalization-group transformations. In Sec.~\ref{sec:conj}, the method to find the conjugation in configuration space is explained. In Sec.~\ref{sec:num}, we provide the numerical results for $d=2$ and $d=3$. 

\section{Renormalization group transformation.}
\label{sec:renorm}

Here the basic idea of renormalization is to construct a canonical change of coordinates iteratively, such that when expressed in the new coordinates, the Hamiltonian has obviously an invariant torus of the chosen frequency. In order to define this change of coordinates, there is a need to enlarge the set of Hamiltonians~\eqref{eqn:Hamiltonian} to the following family of analytic Hamiltonians,
\begin{equation} \label{hamil-fam}
H({\bf A}, {\bm \varphi}) = {\bm \omega} \cdot {\bf A} + \frac{1}{2} \left({\bm \Omega} \cdot {\bf A} \right)^2 +  \sum_{j=0}^\infty f^{(j)}({\bm \varphi})\left({\bm\Omega}\cdot {\bf A}\right)^j.
\end{equation}
Without loss of generality, we can always assume that $\langle f^{(0)}\rangle = 0$, where $\langle g \rangle $ denotes the average of the function $g$ over ${\mathbb T}^d$.
Below, we perform canonical transformations which remain inside this family of Hamiltonians. 

The renormalization transformation we consider follows a scheme proposed in \cite{Koch-99}. This transformation has been used to study the breakup of invariant tori~\cite{Cha-Jau-02, Cha-Las-Ben-Jau-01}. The renormalization transformation, $\R$, is defined for a fixed frequency vector ${\bm \omega} \in \real^d $ and combines a process of rescaling and elimination.  For the sake of completeness, we briefly describe the construction of the renormalization-group transformation.

\subsection{Non-resonant modes }
The cornerstone of the renormalization and its two steps, rescaling and elimination, relies on the distinction between resonant and non-resonant modes. Essentially the resonant (or more precisely near-resonant) modes are the Fourier modes leading to the small denominators in the KAM theory, i.e., the modes $\bm\nu$ of the perturbation where ${\bm\omega}\cdot{\bm\nu}$ are small in amplitude. In other words, the resonant modes are the Fourier modes of the perturbation that influence the motion for a relatively long time. 
We define the set of \textit{non-resonant} modes as follows~\cite{Koch-99,Aba-Koc-Witt-98}
\begin{equation}
\label{nonresonant modes}
    I^- :=\left\{ ({\bm\nu},j) \in I : |{\bm\omega}\cdot {\bm \nu}| > \sigma |{\bm\nu}| + j\kappa \right\}, 
\end{equation}
where $I = \integer^d\times\integer^+$; and $\sigma$, $\kappa$ are fixed constants. The set of resonant modes are defined as $I^+:=I \backslash I^-$. These sets define projection operators $\mathbb{I}^-$ and $\mathbb{I}^+$ as
\begin{eqnarray*}
&& \mathbb{I}^- H = \sum_{({\bm\nu},j)\in I^-} f^{(j)}_{\bm\nu} {\rm e}^{i{\bm\nu}\cdot {\bm\varphi}}({\bm\Omega}\cdot{\bf A})^j, \\
&& \mathbb{I}^+ H = {\bm \omega} \cdot {\bf A} + \frac{1}{2} \left({\bm \Omega} \cdot {\bf A} \right)^2 + \sum_{({\bm\nu},j)\in I^+} f^{(j)}_{\bm\nu} {\rm e}^{i{\bm\nu}\cdot {\bm\varphi}}({\bm\Omega}\cdot{\bf A})^j. 
\end{eqnarray*}

\subsection{Rescaling}

The main purpose of the rescaling procedure is to first move some resonant modes into the non-resonant region. This corresponds to changing the time-scale of interest, moving to larger and larger time-scales. This is achieved by the following canonical transformation:  
\begin{equation}
\label{N}
    \N({\bf A}, {\bm \varphi}) = (N^*{\bf A}, N^{-1} {\bm \varphi} ),
\end{equation}
where $N^*$ denotes the transpose of $N$. 
This transformation moves the modes ${\bm \nu}$ into $N^*{\bm\nu}$. However, this transformation does not leave invariant the family~\eqref{hamil-fam} since it changes ${\bm\omega}\cdot {\bf A}$ into $N{\bm\omega}\cdot {\bf A}=\theta_1 {\bm\omega}\cdot {\bf A}$. In addition, since the terms in  ${\bm\Omega}\cdot {\bf A}$ are also changed into $N{\bm\Omega}\cdot {\bf A}$, and the map ${\bm\Omega}\mapsto N {\bm\Omega}$ is unbounded, some rescalings in time and in action are needed. Following these requirements, we define the rescaling procedure of the renormalization as  
\begin{equation*}
        H'({\bf A}, {\bm \varphi}) = \mu\lambda H(\lambda^{-1} N^*{\bf A}, N^{-1} {\bm \varphi} ).
\end{equation*}
We verify that this transformation leaves the family of Hamiltonians~\eqref{hamil-fam} invariant provided some proper choices for the parameters $\mu$ and $\lambda$: 
$$
    H'({\bf A}, {\bm\varphi}) = \mu \theta_1 {\bm \omega}\cdot {\bf A} + \frac{\mu \| N{\bm\Omega}\|^2}{2\lambda} \left(\frac{ N{\bm\Omega}}{\|N{\bm\Omega}\|}  \cdot{\bf A}\right)^2 + \mu \lambda\sum_{j=0}^\infty f^{(j)}(N^{-1}{\bm\varphi})\left(\frac{\|N{\bm\Omega} \|}{\lambda}\right)^j\left( \frac{N{\bm\Omega}}{\|N{\bm\Omega}\|}\cdot A  \right)^j.
$$
In the above expression, there are two terms contributing to $(N{\bm\Omega}\cdot {\bf A})^2$. We adjust the rescaling $\lambda$ to combine these two terms and keep its value at $1/2$. 
Choosing $\mu = \theta_1^{-1}$, $\lambda = \mu \|N{\bm\Omega}\|^2(1+2\langle f^{(2)}\rangle)$ and denoting 
$$
f'^{(j)}({\bm\varphi}) := \mu\lambda \left(  \frac{\|N{\bm\Omega} \|}{\lambda}  \right)^j  f^{(j)}(N^{-1}{\bm\varphi}), 
$$
we have 
 \begin{equation} \label{second-step}
     H'({\bf A}, {\bm \varphi}) = {\bm \omega}\cdot {\bf A} + \frac{1}{2}({\bm \Omega}'\cdot {\bf A})^2 +  \sum_{j=0}^\infty f'^{(j)}({\bm \varphi})\left({\bm\Omega}'\cdot {\bf A}\right)^j,
 \end{equation}
 where $\langle f'^{(2)}\rangle =0$ and 
 $$
 {\bm\Omega}' = \frac{N{\bm\Omega}}{\|N{\bm\Omega}\|}.
 $$
 \begin{remark}
 As it is pointed out in \cite{Koch-99}, the map $H\mapsto H'$ defined in \eqref{second-step} is not a dynamical system of a space of analytic Hamiltonians to itself. More precisely, defining $\mathcal{A}_\rho$ as the space of analytic Hamiltonians  on $\mathcal{D}_\rho := \{\|A\|< \rho\} \times \{|\mathrm{Im}({\bm\varphi})| < \rho \}$, the map $H \mapsto H'$ is not a dynamical system on any space $\mathcal{A}_\rho$. The latter is due to the fact that the domain $\mathcal{D}_\rho$ is not left invariant by $\N$ in \eqref{N}, in fact $N^{-1}$ expands in the direction of ${\bm \omega}$ which produces a loss of analyticity in the variable ${\bm\varphi}$. To \emph{avoid} this loss of analyticity, the aim is to completely eliminate the non-resonant part of the Hamiltonian by means of a canonical change of coordinates.

 \end{remark}

 \subsection{Elimination}
 
 The second step of the renormalization transformation is a canonical change of coordinates, $\mathcal{U}_{H}$, which eliminates the non-resonant modes of $H$, that is,
     $$\mathbb{I}^-( H\circ \mathcal{U}_{H}) = 0.$$
The construction of this canonical transformation follows a KAM-type procedure. The idea is to construct recursively a sequence of   Hamiltonians $H_k$, with $H_0 = H$, such that the limit $H_\infty$ contains only resonant modes. Each step of the procedure is done by applying a canonical change of coordinates, $H_{k+1} = H_k\circ\U_k$, such that the order of the non-resonant modes of $H_{k+1}$ is $\e_k^2$, i.e.,
\begin{equation}
    \mathbb{I}^-(H_k\circ \mathcal{U}_k)  = {\mathcal O}(\varepsilon_0^{2^{k+1}}),
\end{equation}
where $\e_k$ denotes the order of the non-resonant modes of $H_k$ and $\mathbb{I}^-H_0 = {\mathcal O}(\e_0)$. When this procedure converges it defines a  canonical transformation  \begin{equation}\label{UH}
    \mathcal{U}_{H} = \mathcal{U}_0\circ \mathcal{U}_1\circ \cdots \circ \mathcal{U}_k\circ \cdots,
\end{equation}
such that $\mathbb{I}^-(H \circ \U_{H}) = 0$ . 
In what follows, we describe the construction of the transformation $\mathcal{U}_k$ for one step of this process. The canonical transformations we use are canonical Lie transforms (for a review see ~\cite{Car-81}). These transformations are obtained from a generating function $S({\bf A}, {\bm \varphi})$:
$$
({\bf A}', {\bm \varphi}') = {\rm e}^{-{\mathcal L}_S} ({\bf A}, {\bm \varphi}),
$$
where ${\mathcal L}_S$ is the Liouville operator generated by $S$ acting on $F({\bf A}, {\bm \varphi})$ as ${\mathcal L}_S F=\{F, S\}$, and $\{\cdot,\cdot\}$ is the Poisson bracket
$$ \{F,G\} = \frac{\partial F}{\partial {\bm \varphi}}\cdot \frac{\partial G}{\partial {\bf A}} - \frac{\partial F}{\partial {\bf A}}\cdot \frac{\partial G}{\partial {\bm \varphi}}.$$
These transformations act on the Hamiltonian as
\begin{equation}
\label{eqn:exp}
H'={\rm e}^{{\mathcal L}_S} H.
\end{equation}
This transformation can be seen as a time-1 map of a continuous Hamiltonian flow generated by the Hamiltonian $-S$. 
The generating function is chosen such that it eliminates the order $\e_k$ of $H_k$. Following \cite{Cha-Jau-02}, the family of generating functions is chosen to be 
\begin{equation}\label{S}
    S({\bf A}, {\bm \varphi}) = i \sum_{j=0} ^\infty Y^{j}({\bm \varphi}) ({\bm\Omega} \cdot {\bf A} )^j + a{\bm \Omega}\cdot {\bm \varphi} .
\end{equation}
One way to compute the Hamiltonian in~\eqref{eqn:exp} is by using its expansion
\begin{equation}\label{u1-def}
    H' = \sum_{k=0}^\infty \frac{{\mathcal L}_S^k H }{k!},
\end{equation}
and the recursion ${\mathcal L}_S^k H = {\mathcal L}_S({\mathcal L}_S^{k-1}H)$ (for explicit formulas, see~\cite{Cha-Jau-02}). 

This representation of the change or coordinates has some drawbacks both analytically and numerically. From the analytical point of view, the operator ${\rm e}^{\epsilon {\mathcal L}_S}$ might have some singularities, depending on ${\mathcal L}_S$,  which can yield a radius of convergence, $r$, such that $0< \epsilon<r<1$.

From a numerical point of view, it is known that even when taking the exponential of a matrix or an operator, the computation of the exponential is a delicate procedure, see for example \cite{Mol-Loa-03}. There is no ideal procedure and it should be adapted to the properties of the matrix or the operator whose exponential has to be computed. However, what is known is that considering the series expansion for the numerical computation of the exponential as in \eqref{u1-def}, although practical, is one of the worst way of doing it. Among the several ways to compute numerically an exponential we have chosen the following two ways: 
\\
1) by computing directly the series ${\rm e}^{\epsilon {\mathcal L}_S}H$  with $\epsilon =1$ as in~\eqref{u1-def}, \\
2) by computing $\left({\rm e}^{\frac{1}{m}{\mathcal L}_S}\right)^m H$ using an adaptive method. Assuming we are able to compute Exp$(\epsilon) H := {\rm e}^{\epsilon{\mathcal L}_S} H$, e.g., using a power series as in~\eqref{u1-def}, the adaptive method can be summarized in the following algorithm:
\begin{verbatim}
 Exp_adaptive(epsilon_0)H:
    step = epsilon_0
    If step < min:
        return Exp(step)H
    res1 = Exp(step)H
    res2 = Exp(0.5*step)Exp(0.5*step)H
    if |res1 -res2| < abstol + reltol|res1|:
        return 0.75*res1 + 0.25*res2
    else:
        return Exp_adaptive(0.5*step)Exp_adaptive(0.5*step)H
\end{verbatim}

The first method, referred to as the time-1 method below, has been shown to be effective to approximate the non-trivial fixed point of the renormalization operator $\R$ for $d=2$ \cite{Cha-Jau-02,Cha-Jau-98}.
We have used a second way to compute the exponential of the Liouville operator as a way to improve the first one, and highlight the practical limitations of the first method. Note that one expects  ${\rm e}^{\frac{1}{m} {\mathcal L}_S}$  to have a larger radius of convergence and to be  numerically more stable. 

For both methods, the generating function will be identical. Its aim is to eliminate the order $\e$ of the non-resonant modes of $H$. 
Writing $H$ in the form $H = H_0 +  V$ with 
\begin{equation}
    H_0({\bf A}, {\bm \varphi}) = {\bm \omega} \cdot {\bf A } +\frac{1}{2}({\bm \Omega} \cdot {\bf A})^2 \quad\mbox{and}\quad V({\bf A}, {\bm\varphi}) = \sum_{j=0}^\infty f^{(j)}({\bm\varphi}) ({\bm \Omega}\cdot {\bf A})^j ,
\end{equation}
and assuming ${\mathbb I}^- V\sim {\mathcal O}(\e)$, $S\sim {\mathcal O}(\e)$, if one computes \eqref{u1-def} one has that the order $\e$ of the non-resonant modes of $H'$ is given by the term ${\mathbb I}^-\left( V + \{S, H_0\}\right)$. Therefore the function $S$ is determined by the equation 
$$
    \mathbb{I}^- V + \mathbb{I}^-\{S, H_0\} = 0.
$$
This equation is solved in Fourier space, which yields 
\begin{eqnarray*}
&& Y^{(0)}  = \sum_{{\bm\nu}\in I^-} \frac{f_{{\bm\nu}}^{(0)}}  {{\bm \omega}\cdot {\bm\nu} } {\rm e}^{i{\bm\nu}\cdot {\bm\varphi}},\\
&& Y^{(j)}  = \sum_{{\bm\nu}\in I^-} \frac{1}{{\bm \omega}\cdot {\bm\nu} } (f_{{\bm\nu}}^{(j)} - 2\langle f^{(2)}\rangle {\bm\Omega}\cdot{\bm\nu} Y_{\bm\nu}^{(j-1)} )   {\rm e}^{i{\bm\nu}\cdot {\bm\varphi}}.
\end{eqnarray*}
The constant $a$ in \eqref{S} eliminates the mean value of the linear term in the variable ${\bm \Omega}\cdot {\bf A}$, so that the invariant torus under consideration is located around ${\bm\Omega}\cdot {\bf A} = 0$ in the new coordinates. Therefore, $a$ must satisfy 
$$
    a = - \frac{\langle f^{(1)} \rangle }{2 \|{\bm\Omega}\|^2 \langle f^{(2)}\rangle }. 
$$
The renormalization-group transformation is the defined as \begin{equation}\label{def:ren-trans}
    \R( H ) = H \circ \mathcal{U}_{H}\circ {\mathcal N}.
\end{equation}

\begin{remark}[Trivial fixed point]
Consider the unperturbed Hamiltonian $H^*_0({\bf A}, {\bm \varphi}) ={\bm \omega} \cdot {\bf A} + \frac{1}{2} \left({\bm \Omega}^* \cdot {\bf A} \right)^2  $, with ${\bm \Omega}^*$ a unit eigenvector of $N$ different from ${\bm \omega}$, say $N{\bm\Omega} = \theta_2{\bm\Omega}$. Since $H_0^*$ only contains resonant modes, one can check that $\R(H^*_0) = H^*_0 $. That is,  $H^*_0$ is a fixed point with a scaling in the actions $\lambda = \theta_1^{-1}\theta_2^2$ satisfying $|\lambda| > 1$; meaning that the renormalization focuses on smaller and smaller regions in the actions around ${\bf A}={\bf 0}$.
\end{remark}

We implement numerically the transformation $\R$ defined in \eqref{def:ren-trans} for Hamiltonians~\eqref{hamil-fam}. The approximations we perform in the numerical implementation are of two types: a truncation of the Fourier series of the functions $f^{(j)}$ as follows 
$$
    f^{(j)}({\bm\varphi}) = \sum_{|{\bm\nu}|_\infty \leq L} f_{\bm\nu} ^{(j)} {\rm e}^{i{\bm\nu}\cdot{\bm\varphi}},
$$  
where $|{\bm\nu}|_\infty = \max_i |\nu_i|$, and a truncation in the power series of the actions by neglecting the terms of order $ {\mathcal O}\left( ({\bm\Omega}\cdot {\bf A})^{J+1}     \right)$, which amounts to $(2L+1)^d$ Fourier coefficients for each scalar function $f^{(j)}$. In addition, the action on ${\bm\Omega}$ is a $d-1$ map, independent of the action on the Fourier coefficients. 
This means that we approximate the renormalization map~\eqref{def:ren-trans} by a $(J+1)(2L+1)^d + d-1$-dimensional map. The set of non-resonant modes, \eqref{nonresonant modes}, is defined by fixing the parameters $\sigma = 0.6$ and $\kappa = 0.1$. The codes written in Python 3 are available open source at \texttt{github.com/apbustamante/Renorm}.  All the codes are written in Python 3 using \texttt{NumPy}~\cite{2020NumPy-Array} and \texttt{SciPy}~\cite{2020SciPy-NMeth}.

Given a Hamiltonian $H$ of the form~\eqref{eqn:Hamiltonian}, the main assumption is that if the successive actions of the renormalization operator $\R^n(H)$ -or more precisely its approximate map- on this Hamiltonian converges to a Hamiltonian of the form $H_0({\bf A}) = {\bm\omega}\cdot{\bf A} + \frac{1}{2}({\bm\Omega}_\infty \cdot{\bf A})^2 $, then $H$ has a smooth invariant torus with frequency $\bm\omega$. If the action of the renormalization map on $H$ diverges, i.e., $\R^n(H)\longrightarrow \infty$ as $n\rightarrow \infty$ , then $H$ does not have this invariant torus. 

One of the objectives of this article is to inspect this assumption. In order to check if $H$ has an invariant torus or not, we compare the renormalization results with a different method. The one we choose follows the KAM theory in configuration space~\cite{Sal-Zeh-89,Cel-Chi-88}. 

\section{Conjugation in configuration space}
\label{sec:conj}

In order for Hamiltonian~\eqref{eqn:Hamiltonian} to have an invariant torus with frequency vector ${\bm\omega}$, we are looking for a conjugation of the type
\begin{eqnarray}
&& {\bf A} = \overline{\bf A}(z,{\bm\psi}),\\
&& {\bm \varphi} = {\bm\psi}+{\bm\Omega} h({\bm\psi}), 
\end{eqnarray}
where the flow in the new coordinates $(z,{\bm\psi})$ is linear, i.e., $\dot{z}=0$ and $\dot{\bm\psi}={\bm\omega}$. The equation of motion for Hamiltonian~\eqref{eqn:Hamiltonian} leads to the following equation for the function $h:{\mathbb T}^d \to {\mathbb R}$ (see also Refs.~\cite{Sal-Zeh-89,Cel-Chi-88}:
\begin{equation}
\label{eqn:confeq}
    \left({\bm\omega}\cdot \frac{\partial}{\partial {\bm\psi}}\right)^2 h({\bm\psi})+{\bm\Omega}\cdot \frac{\partial V}{\partial {\bm\varphi}}\left({\bm\psi}+{\bm\Omega} h({\bm\psi}) \right)=0.
\end{equation}
\begin{remark}[Gauge symmetry]
 If $h({\bm \psi})$ is a solution of Eq.~\eqref{eqn:confeq}, then $\hat{h}({\bm \psi})=h({\bm \psi} +\eta {\bm\Omega})+\eta$ is also a solution for all $\eta$. We consider solution of Eq.~\eqref{eqn:confeq} with $\langle h\rangle =0$. 
\end{remark}
We use a Newton method as developed in Refs.~\cite{Su-del-12,Bla-del-13} to solve Eq.~\eqref{eqn:confeq}. More explicitly, we consider the following equation
$$
\left({\bm\omega}\cdot \frac{\partial}{\partial {\bm\psi}}\right)^2 h({\bm\psi})+{\bm\Omega}\cdot \frac{\partial V}{\partial {\bm\varphi}}\left({\bm\psi}+{\bm\Omega} h({\bm\psi}) \right)+\lambda=0.
$$
We assume that after $n$ steps of the Newton method, we have approximate solutions $h_n$ and $\lambda_n$ such that
\begin{equation}
\label{eqn:h0}
  \varepsilon({\bm\psi})=\left({\bm\omega}\cdot \frac{\partial}{\partial {\bm\psi}}\right)^2 h_n({\bm\psi})+{\bm\Omega}\cdot \frac{\partial V}{\partial {\bm\varphi}}\left({\bm\psi}+{\bm\Omega} h_n({\bm\psi}) \right)+\lambda_n,  
\end{equation}
is small. We are looking for a refined solution $h_{n+1}=h_n+\Delta$ and $\lambda_{n+1}=\lambda_n+\delta$ such that $h_{n+1}$ is closer to a true solution of Eq.~\eqref{eqn:confeq}. The increments are now solution of 
$$
\left({\bm\omega}\cdot \frac{\partial}{\partial {\bm\psi}}\right)^2 \Delta+{\bm\Omega}\cdot \frac{\partial V}{\partial {\bm\varphi}}\left({\bm\psi}+{\bm\Omega} (h_n+\Delta) \right)-{\bm\Omega}\cdot \frac{\partial V}{\partial {\bm\varphi}}\left({\bm\psi}+{\bm\Omega} h_n \right) +\delta = -\varepsilon.
$$
If $\Delta$ is sufficiently small, a good approximation for $\Delta$ can be defined by
$$
\left({\bm\omega}\cdot \frac{\partial}{\partial {\bm\psi}}\right)^2 \Delta+\Delta \left({\bm\Omega}\cdot \frac{\partial }{\partial {\bm\varphi}}\right)^2 V\left({\bm\psi}+{\bm\Omega} h_n \right)  +\delta = -\varepsilon,
$$
neglecting the second-order terms in $\Delta$.
The second-order derivative can be computed by differentiating Eq.~\eqref{eqn:h0}:
$$
{\bm\Omega}\cdot \frac{\partial \varepsilon}{\partial {\bm\psi}}=\left({\bm\omega}\cdot \frac{\partial}{\partial {\bm\psi}}\right)^2 l({\bm\psi})+l({\bm\psi }) \left( {\bm\Omega}\cdot \frac{\partial }{\partial {\bm\varphi}}\right) ^2 V\left({\bm\psi}+{\bm\Omega} h_n({\bm\psi}) \right),
$$
where $l({\bm\psi}) = 1+ {\bm\Omega}\cdot {\partial h_n}/{\partial {\bm\psi}}$.
Neglecting terms of order $\varepsilon \Delta$, the equation for $\Delta$ becomes
\begin{equation} 
\label{eqn:lDel}
l({\bm\psi}) \left({\bm\omega}\cdot \frac{\partial}{\partial {\bm\psi}}\right)^2 \Delta -\Delta \left({\bm\omega}\cdot \frac{\partial}{\partial {\bm\psi}}\right)^2 l({\bm\psi})=-l({\bm\psi})\left(\delta+\varepsilon({\bm\psi})\right).
\end{equation}
We solve Eq.~\eqref{eqn:lDel} using two cohomological equations for the auxiliary functions $W$ and $\beta$ (chosen with zero meanvalue):
\begin{eqnarray}
&& {\bm\omega}\cdot \frac{\partial W}{\partial {\bm\psi}}=l({\bm\psi})\left(\delta+\varepsilon({\bm\psi})\right),\label{eqn:coho1}\\
&& {\bm\omega}\cdot \frac{\partial \beta}{\partial {\bm\psi}}=-\frac{W+W_0}{l({\bm\psi})^2},\label{eqn:coho2}
\end{eqnarray}
where 
\begin{eqnarray}
&& W_0=-\frac{\langle W/l^2\rangle}{\langle 1/l^2\rangle},\\
&& \delta = -\langle l \varepsilon\rangle,
\end{eqnarray}
in order to ensure the existence of a solution for Eqs.~\eqref{eqn:coho1}-\eqref{eqn:coho2}. Furthermore since the solutions of Eq.~\eqref{eqn:lDel} are defined up to a constant multiplied by $l$, we impose $\langle \Delta \rangle =0$ to fix the gauge. Therefore the solution of Eq.~\eqref{eqn:lDel} for $\Delta$ is given by
$$
\Delta = \beta l- l \langle \beta l \rangle. 
$$
Numerically, Eqs.~\eqref{eqn:coho1}-\eqref{eqn:coho2} are solved using Fourier transforms (see also Refs.~\cite{Call-del-09,Call-del-10}). The new solution is defined by
$$
h_{n+1}=h_n+\Delta,
$$
together with $\lambda_{n+1}=\lambda_n+\delta$.  
Provided that the potential is sufficiently small, we initiate the Newton method with the following initial guess: 
\begin{equation}
\label{eqn:hinit}
h_0=-\left({\bm\omega}\cdot \frac{\partial}{\partial {\bm\varphi}}\right)^{-2} {\bm\Omega}\cdot \frac{\partial V}{\partial {\bm\varphi}},
\end{equation}
and $\lambda_0=0$.
Other strategies to design better initial guesses might be more suitable, involving, e.g., a continuation method from an integrable case as in Refs.~\cite{Call-del-10,Bla-del-13} or an expansion in the small parameters as in Ref.~\cite{Cel-Chi-88}. We have chosen the initial guess~\eqref{eqn:hinit} as a common reference point for counting the number of iterations of the Newton method to reach a good accuracy of the solution of Eq.~\eqref{eqn:confeq}. 

If the Newton iterations converge, i.e., if there exists $n_*$ such that $\Vert h_{n_*+1}-h_{n_*}\Vert\leq \eta$ where $\eta$ is a small threshold parameter, it can be proved that the Hamiltonian system has an invariant torus with frequency vector ${\bm \omega}$ (see Refs.~\cite{Sal-Zeh-89,Cel-Chi-88}). 

In the numerical implementation, we project $h({\bm \psi})$ in Fourier space with $L^d$ Fourier modes, and use extensively fast Fourier transforms. We use the monitoring of the tail in Fourier series in order to adjust the value of $L$ (for more details, see Ref.~\cite{Har-Can-Fig-16}). We also remove the Fourier modes with amplitudes smaller than a certain threshold to remove some numerical instability (see Ref.~\cite{Bla-del-13}). The codes written in Python 3 are available open source at \texttt{github.com/cchandre/ConfKAM}. All the codes are written in Python 3 using \texttt{NumPy}~\cite{2020NumPy-Array}.

\section{Numerical results}
\label{sec:num}

In appearance, the conjugation in configuration space is much simpler to implement since it amounts to the determination of a single scale function of $d$ angles, and that its numerical implementation does not dependent on the number-theoretic properties of the frequency vector. As for the renormalization, it involves the determination of $J+1$ scalar functions of $d$ angles, and its implementation is tailored to the frequency vector through the choice of the matrix $N$ and the set of resonant/non-resonant modes. It is therefore more difficult to adapt it to a generic frequency vector when compared with the conjugation in configuration space. 
The question regarding which method is numerically more efficient boils down to how many Fourier modes are needed to represent the scalar functions (up to some given accuracy). Given that there is a specific treatment of the resonant modes, we expect that the renormalization will need fewer modes. However this depends on the number-theoretic properties of the frequency vector.  
In order to bring some elements of answer, we apply the two methods, the renormalization and the conjugation in configuration space, for $d=2$ and $d=3$ for some rather simpler frequency vectors, one related to the golden mean in 2D, and one related to the spiral mean in 3D. The conclusions drawn in this section cannot be generalized to other frequency vectors, but these examples shed some light on the advantages and complementarity of both methods. 

\subsection{For $d=2$}

We consider the following Hamiltonian
\begin{equation}
\label{eqn:numHam2D}
    H({\bf A}, {\bm \varphi})={\bm \omega}\cdot {\bf A} +\frac{1}{2}({\bm \Omega}\cdot {\bf A})^2+\mu_1 \cos \varphi_1 +\mu_2 \cos (\varphi_1+\varphi_2),
\end{equation}
with ${\bm \Omega}=(1, 0)$. We consider the frequency vector ${\bm \omega}=(\omega, -1)$ with $\omega=(\sqrt{5}-1)/2$, which is an eigenvector of $N=\left( \begin{array}{cc}
  1   & 1 \\
    1 & 0
\end{array}\right)$ with eigenvalue $-\omega$. 
The objective is to determine the set of parameters $(\mu_1,\mu_2)$ for which Hamiltonian~\eqref{eqn:Hamiltonian} has an invariant torus with the given frequency vector. 

\begin{figure}
    \centering
\includegraphics[width = 15.0truecm]{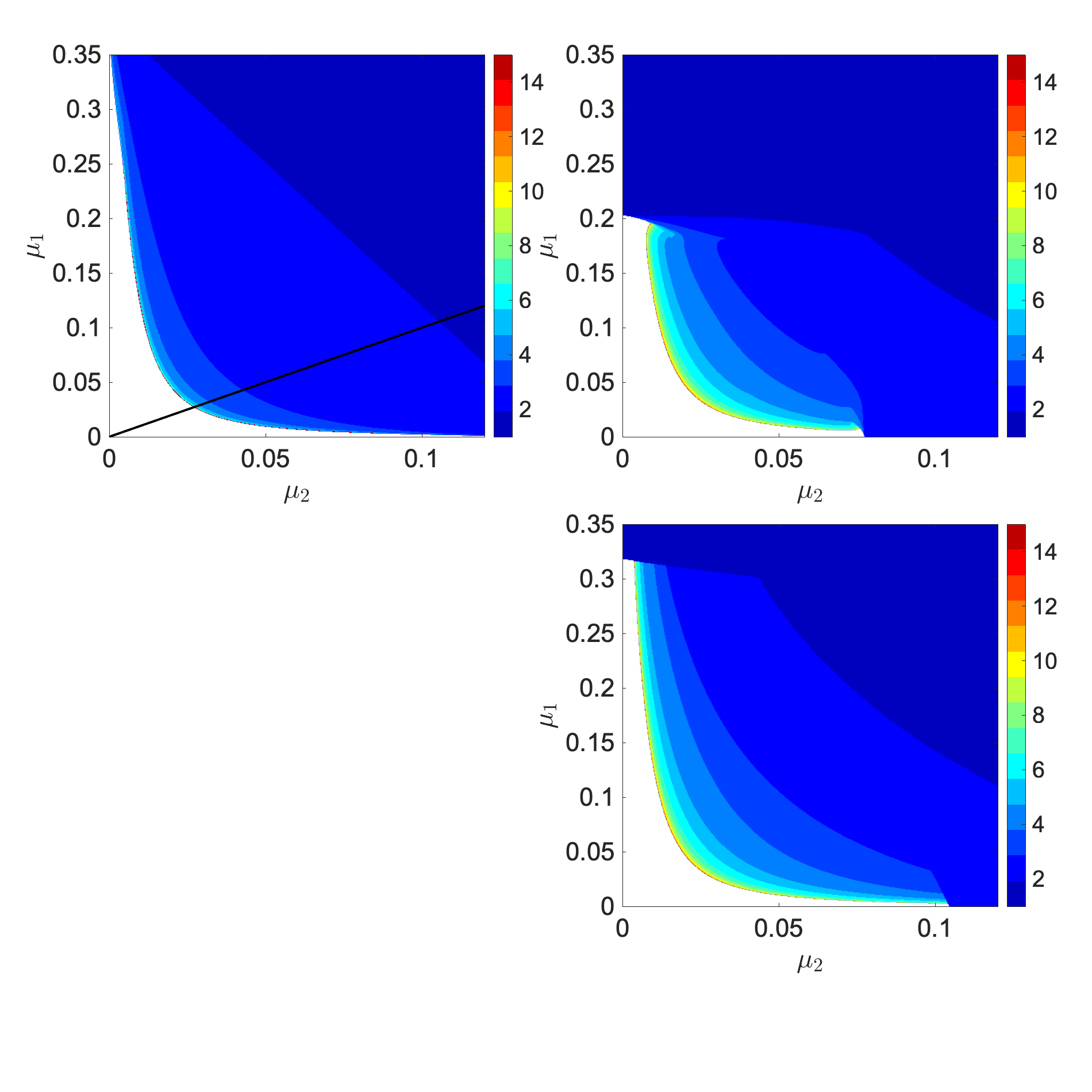}
\caption{Critical surfaces for the invariant torus of Hamiltonian~\eqref{eqn:numHam2D} with frequency $\omega=(\sqrt{5}-1)/2$. Top left: conjugation in configuration space (with $L=2^{10}$). Top right: renormalization with time-1 Lie transforms ($L=J=5$). Bottom right: renormalization with adaptive time-step Lie transform ($L=J=5$). The white region represents the domain of convergence of the iterates of the methods, and the color-scaled region, the number of iterates before the norm exceeds a threshold for divergence ($10^2$ for the conjugation method and $10^{10}$ for the renormalization). } 
\label{alpha1}
\end{figure}

In Fig.~\ref{alpha1}, the domain of convergence of the iterates of the renormalization map and the conjugation method is represented for Hamiltonians~\eqref{eqn:numHam2D}. If $(\mu_1,\mu_2)$ is in the white region, this means that the renormalization or the conjugation method was able to find numerically the conjugation to a trivial system which has an invariant torus of the chosen frequency. The white region is then expected to be the region where Hamiltonian~\eqref{eqn:numHam2D} has an invariant torus with the chosen frequency. As expected from KAM theory, this region contains the region around $\mu_1=\mu_2=0$. Moreover, since Hamiltonian~\eqref{eqn:numHam2D} is integrable if $\mu_1=0$ or $\mu_2=0$, the white regions should include these lines in parameter space and a region around them. 

Using a color scale, we represent the number of iterates necessary for the method to exceed a given threshold ($10^8$ for the conjugation method and $10^{10}$ for the renormalization). This color scheme highlights better the critical surface, i.e., the surface where the methods do neither converge to a trivial system and neither diverge. As expected, it shows that the two methods need more iterates close to the critical surface. We also notice that it takes about the same number of iterations to exceed the divergence threshold, and this number is rather low, so the divergence above the critical surface is rather fast for both methods. 

We notice that the conjugation in configuration space has a wider domain of convergence compared to the renormalization. In fact, the critical surface obtained by the renormalization is distorted by a domain in which divergence occurs in one or two iterations (dark blue region in Fig.~\ref{alpha1}, top right panel). This corresponds to a domain where one of the canonical transformations define by a time-1 Lie transform does not converge. In order to remedy this problem, we have used a variable time-step Lie transform in Fig.~\ref{alpha1} (lower right panel). Indeed we check that the critical surface corresponds approximately to the one obtained using the conjugation in configuration space (except at the extremities where the canonical transformations diverge again). The adaptive step size method for the computation of the Lie transforms is able to increase the domain of convergence of the renormalization map. In a wide region of the parameter space $(\mu_1,\mu_2)$, we observe some good quantitative agreement, which indicates that both methods, if converging, converge in the domain of existence of the invariant torus of the chosen frequency. As expected, the observed quantitative agreement depends on the parameters of the models.

The computation of the Lie transforms is the bottleneck for the convergence of the renormalization operator. The bottleneck for the conjugation in configuration space is the number of Fourier modes necessary to represent $h({\bm\psi})$ in order to accurately represent invariant tori close to the breakup. 
More quantitatively, for $\mu_1=\mu_2=\epsilon$, the critical threshold obtained by renormalization is $\epsilon_c\approx 0.027590$ (obtained with a relatively low number of Fourier modes $L=5$), whereas $\epsilon_c\approx 0.027509$ with the conjugation in configuration space with $L=2^{13}$ (and $\epsilon_c\approx 0.026909$ for $L=2^{10}$ for comparison). The critical value obtained using the conjugation method is always smaller than the critical threshold for the breakup of the invariant torus under consideration. For this case in 2D, the value is rather close to the critical value. We notice that it was proven in Ref.~\cite{Cel-Gio-Loc-00} that an invariant torus exists for $\epsilon < 0.025375$. The correct value obtained by other methods like Greene's residue criterion is $\epsilon_c\approx 0.027590$ (see Ref.~\cite{Cha-Jau-02}). 

\subsection{For $d=3$}

We consider the following frequency vector ${\bm \omega}=(\sigma, \sigma^2, 1)$ where $\sigma\approx 1.324718$ is the real root of 
$$
\sigma^3 = \sigma+1. 
$$
It is an eigenvector of $N=\left( \begin{array}{ccc}
    0 & 0 & 1 \\
    1 & 0 & 0 \\
    0 & 1 & -1
\end{array}\right)$ with eigenvalue $1/\sigma$.
This frequency vector has been considered in Refs.~\cite{Art-Cas-She-91,Art-Cas-She-92,Cha-Las-Ben-Jau-01,Cha-Jau-98}.
The Hamiltonian family we choose is
\begin{equation}
\label{eqn:numHam3D}
    H({\bf A}, {\bm \varphi})={\bm \omega}\cdot {\bf A} +\frac{1}{2}({\bm \Omega}\cdot {\bf A})^2+\mu_1 \cos \varphi_1+\mu_2 \cos \varphi_2 +\mu_3 \cos \varphi_3.
\end{equation}
In what follows, we fix $\mu_3=0.1$. The vector ${\bm\Omega}$ is chosen as ${\bm\Omega}=(1,1,-1)$.

\begin{figure}[ht]
    \centering
\includegraphics[width = 15.0truecm]{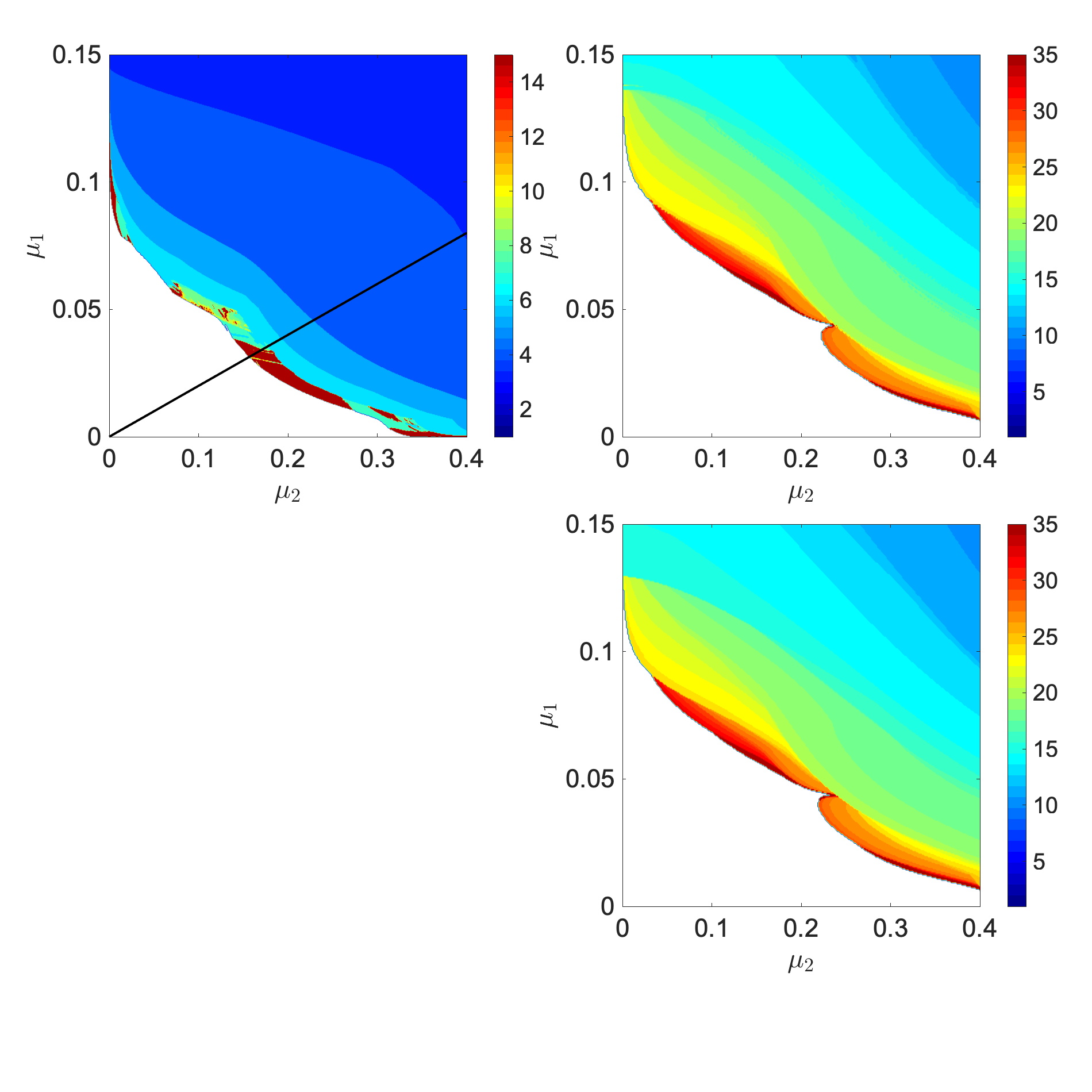}
\caption{Critical surfaces for the invariant torus of Hamiltonian~\eqref{eqn:numHam3D} with frequency vector ${\bm \omega}=(\sigma,\sigma^2,1)$. Top left: conjugation in configuration space (with $L=2^7$). Top right: renormalization with time-1 Lie transforms (with $L=J=5$). Bottom right: renormalization with adaptive time-step Lie transform (also with $L=J=5$). The white region represents the domain of convergence of the iterates of the methods, and the color-scaled region, the number of iterates before the norm exceeds a threshold for divergence ($10^2$ for the conjugation method and $10^{10}$ for the renormalization).} 
\label{3d}
\end{figure}

For the renormalization map, we choose $L=J=5$. For the conjugation in configuration space, we use $2^{21}$ Fourier modes, i.e., $L=2^7$.  In Fig.~\ref{3d}, we represent the domains of convergence and divergence of both methods in the space of parameters $(\mu_1, \mu_2)$. The main striking feature is that these domains do not coincide, even approximately, when comparing the two methods. In particular, the domain of convergence of our implementation of  the method of conjugation in configuration space is significantly smaller than the one of the renormalization. This was also what was observed in 2D, but here the differences are more striking. More quantitatively, we consider a one-parameter family for $\mu_1=\mu_2/5=\epsilon$ and $\mu_3=0.1$ (represented by a continuous black curve in Fig.~\ref{3d}). The critical value given by renormalization is approximately $\epsilon\approx 0.04468$. The method of conjugation fails to provide a solution for $L=2^{9}$ at $\epsilon \gtrsim 0.036353$ (and $\epsilon \gtrsim 0.030226$ for $L=2^{7}$, and  $\epsilon \gtrsim 0.035160$ for $L=2^{8}$ for comparison). 
As expected, the main bottleneck of the conjugation in configuration space, namely the number of Fourier modes, turns out to be a major one in 3D. 
\begin{figure}
    \centering
\includegraphics[width = 15.0truecm]{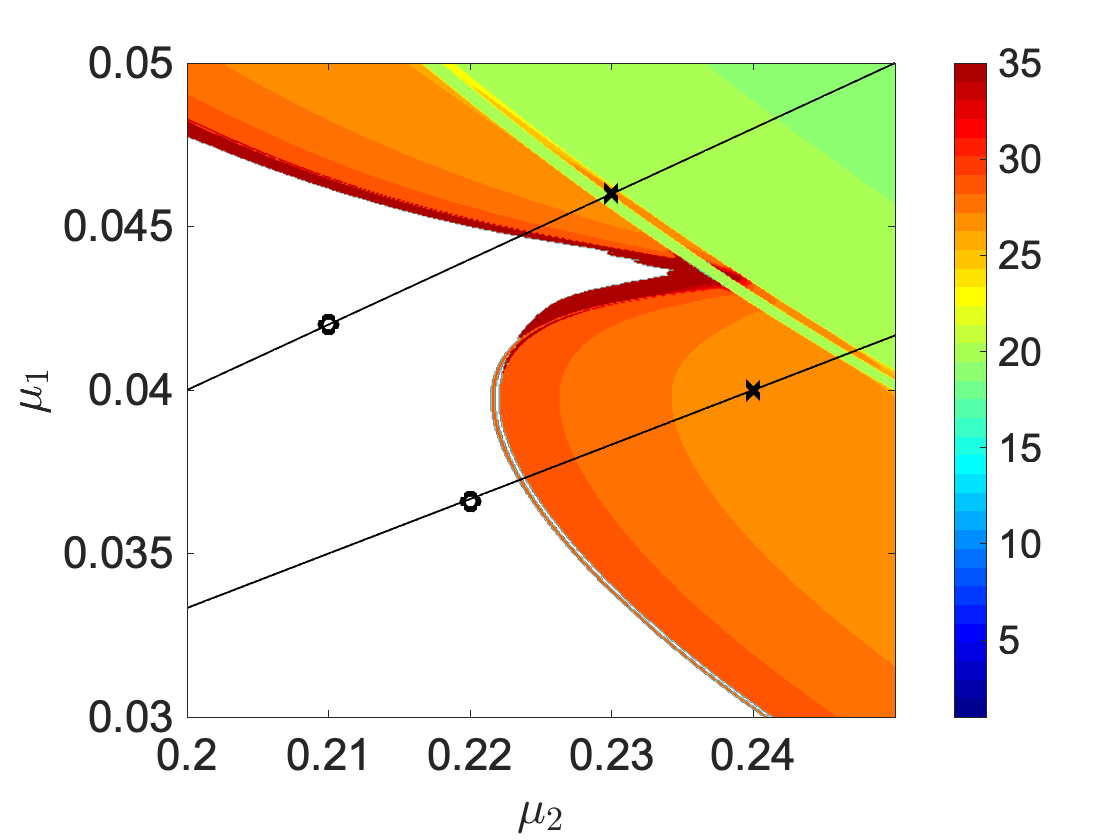}
\caption{Zoom from Fig.~\ref{3d} (upper right panel). The two lines are the two families $\mu_1 = \mu_2/5$ and $\mu_1 = \mu_2/6$. The circles are located at $(\mu_1,\mu_2)=(0.042,0.21)$ and $(\mu_1,\mu_2)=(0.0366,0.22)$. The crosses are located at $(\mu_1,\mu_2)=(0.046,0.23)$ and $(\mu_1,\mu_2)=(0.04,0.24)$. } 
\label{3dzoom}
\end{figure}
Another striking feature is observed when we compare the renormalization results between two and three dimensions: the number of iterations of the renormalization map to diverge is significantly higher for three dimensions than it is for two dimensions. In particular, close to the critical surface (10\% above it), a divergence is observed for 10 iterations in 2D, whereas 30 or more iterations are needed in 3D. It is also worth noticing that the method of conjugation in configuration space diverges fast when it diverges. 
The slow divergence of the renormalization map indicates a rather puzzling renormalization dynamics on the critical surface, which was hinted as a critical non-chaotic strange attractor in Refs.~\cite{Cha-Jau-98,Cha-Jau-Ben-Cell-99}. The critical surface seems to be more sticky for the renormalization dynamics in 3D than it is in 2D. In addition, the critical surface obtained using the renormalization presents some cusps, one being clearly visible around $(\mu_1,\mu_2) \approx (0.044,0.23)$ (see also Fig.~\ref{3dzoom}). Given the high number of iterates necessary to diverge in the vicinity of the cusps, they are not due to a singularity in the definition of canonical Lie transforms. The presence of these cusps is confirmed using variable time-step Lie transforms. Therefore these cusps are not due to the lack of convergence of Lie transforms as it was the case in 2D close to the integrable lines. 

Since the two methods, renormalization and conjugation in configuration space, cannot be compared due to the limits of our implementation of the latter one in 3D, we use another method to assess whether or not the renormalization method converges up to the critical threshold of break-up of the invariant torus under consideration. This method consists in computing rotation numbers using weighted Birkhoff averages.

Given an homeomorphism of the circle $f:\torus \longrightarrow \torus $ and a lift $F:\real \longrightarrow\real$, the rotation number of $f$ is  defined as $$ \rho(f) = \lim_{n\rightarrow \infty} \frac{F^n(\theta) - \theta}{n} ,  $$ for any $\theta\in \torus$.  The computation of $\rho$ using weighted Birkhoff averages  is known as a useful technique to distinguish between chaotic and quasiperiodic dynamics \cite{Mei-San-20, Das-Yor-18}.
Given any function $h:\torus\rightarrow\real$, the Birkhoff average of the lift $F$ is defined as 
$$
    B_S(h)(\theta)= \sum_{n=0}^{S-1}h\circ F^n(\theta),
$$ 
and the weighted Birkhoff average of $F$ as 
$$
    WB_S(h)(\theta) = \frac{1}{C_S}\sum_{n=1}^{S-1}w\left(\frac{n}{S}\right)h\circ F^n(\theta),
$$
where $C_S:= \sum_{n=1}^{S-1}w({n}/{S}) $ and $w$ is the \textit{bump} function 
$$
w(t) = \exp \left(- \frac{1}{t(1-t)} \right).
$$
It is know that $WB_S(h)(\theta)$ converges to $\lim_{n\rightarrow\infty} B_S(h)(\theta)$, when the limit exists, and that the speed of convergence of $WB_S(\theta)$ is super-polynomial, i.e., faster than any powers of $S^{-1}$ \cite{Das-Yor-18}.

Applying a linear change of coordinates $({\bf A}', {\bm \varphi}') = (M{\bm A}, \Tilde{M}^{-1}{\bm \varphi}) $ with
$$
M=\left(\begin{array}{ccc} 1 & 1 & -1\\ 0 & 1 & 0 \\ 1 & 0 & 0\end{array} \right),
$$
Hamiltonian~\eqref{eqn:numHam3D} is mapped into 
\begin{equation}\label{nonautonomous}
    \widetilde{H}({\bf A}, {\bm \varphi})= \frac{1}{2} A_1^2 - A_1 + \nu_1 A_2 + \nu_2 A_3 +\mu_1 \cos( \varphi_1 +\varphi_3) +\mu_2 \cos(\varphi_1 + \varphi_2) +\mu_3 \cos \varphi_1, 
\end{equation}
with $ \nu_1 = \sigma^2 + 1$, $\nu_2 = \sigma+1$. So, the existence of an invariant torus of frequency ${\bm\omega}$ for Hamiltonian~\eqref{eqn:numHam3D} is equivalent to the existence of an invariant torus of frequency $(-1, \nu_1, \nu_2)$ for $\widetilde{H}$. The equations of motion for Hamiltonian $\widetilde{H}$ can be reduced to the following system:
\begin{subequations}
\label{eq-motion}
\begin{align}
   & \frac{{\rm d} \varphi}{{\rm d}t}  = A -1,  \\
   & \frac{{\rm d} A}{{\rm d}t} = \mu_1 \sin( \varphi +\nu_2 t) +\mu_2 \sin(\varphi + \nu_1 t ) +\mu_3 \sin \varphi. 
\end{align}
\end{subequations}
We define the stroboscopic map $F: \torus\times \real \longrightarrow \torus\times \real$ as the time-$2\pi/\nu_2$ map obtained from Eqs.~ \eqref{eq-motion}. We compute the rotation number as a function of the initial conditions $(0, A_0)$. The invariant curve is expected to be at $\rho=-1/\nu_2$. 

\begin{figure}
    \centering
    \includegraphics[width = 18.0truecm]{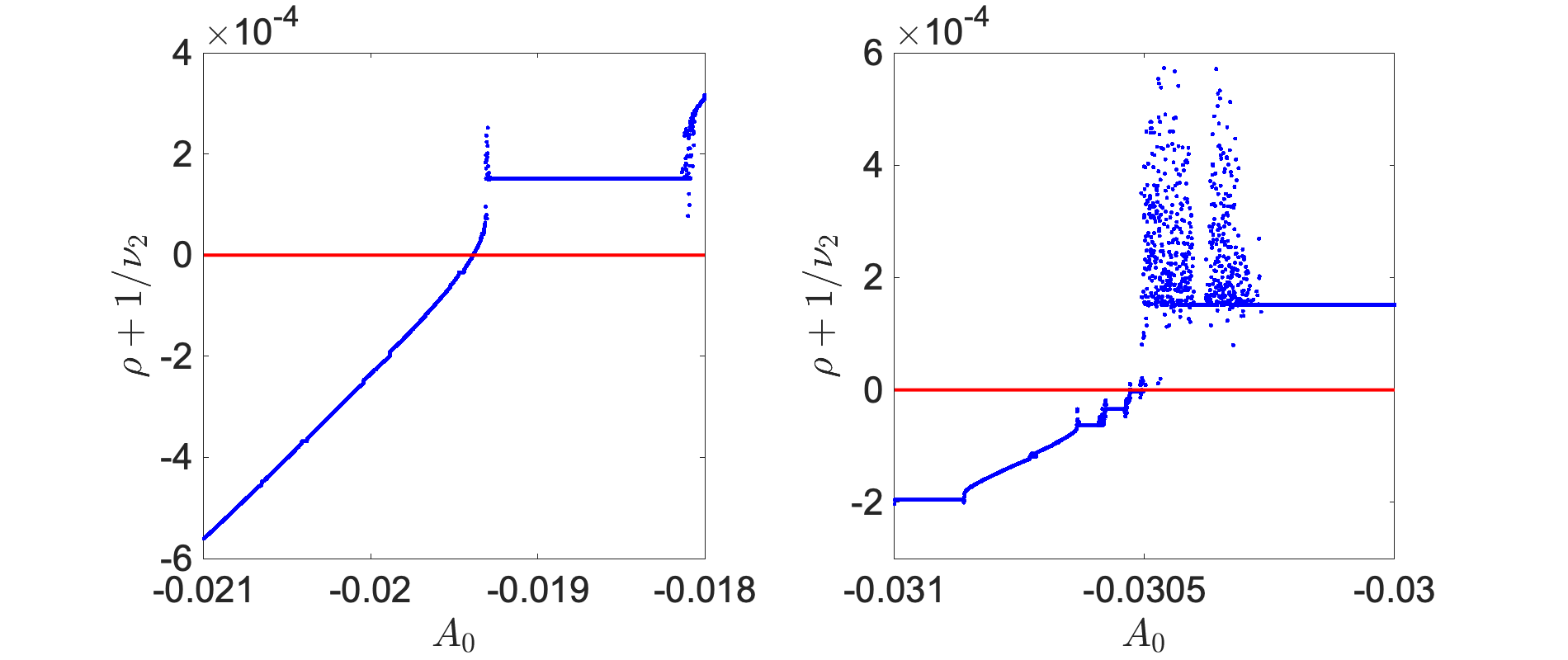}
    \caption{Plot of the rotation number $\rho$ of the stroboscopic map $F$ (see text) as a function of the initial condition $A_0$ for different values of the parameters. Left panel: $(\mu_1,\mu_2)=(0.042,0.21)$ (upper circle on Fig.~\ref{3dzoom}). Right panel: $(\mu_1,\mu_2)=(0.046,0.23)$ (upper cross on Fig.~\ref{3dzoom}). The rotation numbers are computed from orbits of length $4 \times 10^4$ points for $F$.}
    \label{rotnumber1}
\end{figure}

\begin{figure}
    \centering
    \includegraphics[width = 18.0truecm]{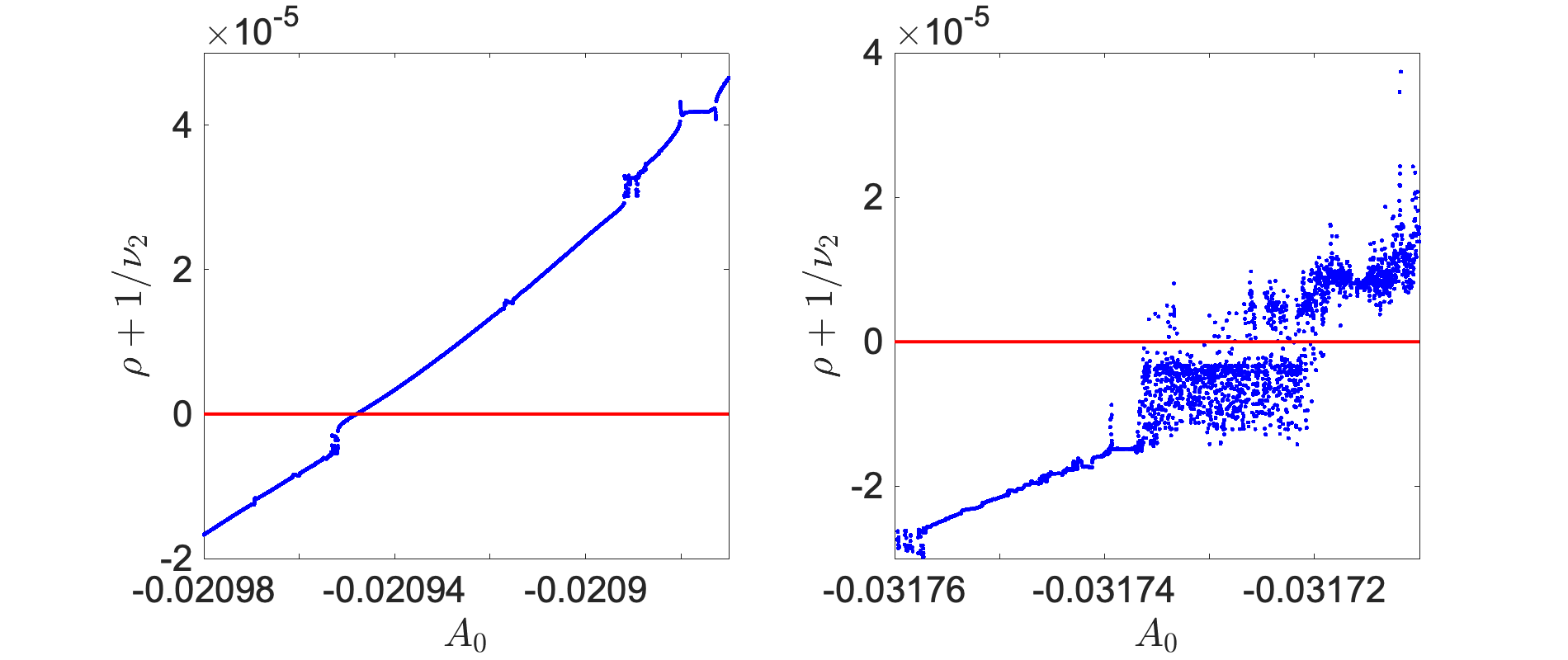}
    \caption{Plot of the rotation number of the stroboscopic map $F$ (see text) as a function of the initial condition $A_0$ for different values of the parameters. Left panel: $(\mu_1,\mu_2)=(0.0366,0.22)$ (lower circle on Fig.~\ref{3dzoom}). Right panel: $(\mu_1,\mu_2)=(0.04,0.24)$ (lower cross on Fig.~\ref{3dzoom}). The rotation numbers are computed from orbits of length $4 \times 10^4$ points for $F$.}
    \label{rotnumber2}
\end{figure}

In Figs.~\ref{rotnumber1} and \ref{rotnumber2} we include the computations of the rotation numbers for some values of  $\mu_1$ and $\mu_2$ close to the \textit{breakup} predicted by the renormalization-group method. The values we consider are represented with circles and crosses on Fig.~\ref{3dzoom}. 
From these computations, we see that when the parameters $(\mu_1,\mu_2)$ are within the domain of convergence of the renormalization map, the corresponding Hamiltonian seems to have an invariant torus from the analysis of the rotation numbers. When the parameters are outside this domain of convergence, the analysis clearly indicates that the Hamiltonian system does not have this invariant torus. The analysis of rotation numbers is conclusive if the system does not have an invariant torus (because of the sensitivity of the frequency with the initial condition), but it is more delicate to conclude the existence of this invariant torus, since it might just be a scale issue, i.e., at a smaller scale than what is shown, the analysis could very well show a sensitivity to initial condition, and hence the absence of invariant torus. 

\section{Conclusions}

Given the specific treatment of the resonant/non-resonant modes, the renormalization method necessitates much fewer modes for each angle than the conjugation in configuration space. When approaching criticality, the number of Fourier modes necessary to describe $h({\bm\psi})$ diverges quickly, illustrating the roughness of the shape of the critical tori~\cite{Har-Can-Fig-16}. As a consequence, the conjugation method is not well suited for the investigation of higher-dimensional tori. Ideally, it would be interesting to combine the advantage of the conjugation method (namely the determination of a single scalar function) and the specific treatment of resonant/non-resonant modes of the renormalization.  

Using the renormalization-group method, we have unveiled the presence of cusps in the critical surface for the breakup of three-dimensional tori. Iterating the renormalization map on the critical surface to derive universal features or scaling relations for the breakup of invariant tori has to take into account these cusps. \\
We notice that some non-smooth features are also visible for the breakdown of analyticity for Frenkel-Kontorova models in quasi-periodic media with two frequencies: We recomputed Fig.~1(A) of Ref.~\cite{Bla-del-13} in a small region around $(\varepsilon_1,\varepsilon_2)\approx (0.01, 0.003)$ (figure not shown), and it showed cusps of a similar nature as the one we unveiled here\footnote{The Python code to produce this figure is available at \texttt{github.com/cchandre/Quasiperiodic\_Frenkel-Kontorova}.}.

\section*{Acknowledgements}
We thank Rafael de la Llave for many comments and suggestions. We also thank the computer support of the School of Mathematics at Georgia Tech.
The project leading to this research has received funding from the European Union’s Horizon 2020 research and innovation program under the Marie Sk\l odowska-Curie Grant Agreement No. 734557. A.P.B. has been partially supported by NSF grant DMS 1800241 and Sloan grant FG-2020-13337.  

\newcommand{\etalchar}[1]{$^{#1}$}

\Addresses
\end{document}